\documentclass[12pt]{amsart}
\usepackage{amsmath}
\usepackage{amsxtra}
\usepackage{amscd}
\usepackage{amsthm}
\usepackage{amsfonts}
\usepackage{amssymb}
\usepackage{eucal}
\usepackage{hyperref}

\def\C{{\mathbb C}}

\def\P{{\mathbb P}}

\def\Z{{\mathbb Z}}

\newtheorem{theorem}{Theorem}[section]

\newtheorem{lemma}[theorem]{Lemma}

\newtheorem{proposition}[theorem]{Proposition}
\newtheorem{corollary}[theorem]{Corollary}

\title{On higher syzygies of ruled surfaces II}
\author{Euisung Park }
\address {Euisung Park : School of Mathematics, Korea Institute for Advanced Study,
207-43 Cheongryangri-dong, Dongdaemun-gu, Seoul 130-722, Republic
of Korea,} \email{puserdos@kias.re.kr}

\begin{document}

\thispagestyle{empty} \maketitle

\setcounter{page}{1}

\begin{abstract}
In this article we we continue the study of property $N_p$ of
irrational ruled surfaces begun in \cite{ES}. Let $X$ be a ruled
surface over a curve of genus $g \geq 1$ with a minimal section
$C_0$ and the numerical invariant $e$.

When $X$ is an elliptic ruled surface with $e = -1$, there is an
elliptic curve $E \subset X$ such that $E \equiv 2C_0 -f$. And we
prove that if $L \in \mbox{Pic}X$ is in the numerical class of
$aC_0 +bf$ and satisfies property $N_p$, then $(C,L|_{C_0})$ and
$(E,L|_E)$ satisfy property $N_p$ and hence $a+b \geq 3+p$ and
$a+2b \geq 3+p$. This gives a proof of the relevant part of
Gallego-Purnaprajna' conjecture in \cite{GP2}.

When $g \geq 2$ and $e \geq 0$ we prove some effective results
about property $N_p$. Let $L \in \mbox{Pic}X$ be a line bundle in
the numerical class of $aC_0 +bf$. Our main result is about the
relation between higher syzygies of $(X,L)$ and those of
$(C,L_{C})$ where $L_C$ is the restriction of $L$ to $C_0$. In
particular, we show the followings: $(1)$ If $e \geq g-2$ and
$b-ae \geq 3g-2$, then $L$ satisfies property $N_p$ if and only if
$b-ae \geq 2g+1+p$. $(2)$ When $C$ is a hyperelliptic curve of
genus $g \geq 2$, $L$ is normally generated if and only if $b-ae
\geq 2g+1$ and normally presented if and only if $b-ae \geq 2g+2$.
Also if $e \geq g-2$, then $L$ satisfies property $N_p$ if and
only if $a \geq 1$ and $b-ae \geq 2g+1+p$.
\end{abstract}

\tableofcontents

\section{Introduction}
\noindent In this article we continue the study of property $N_p$
of irrational ruled surfaces begun in \cite{ES}. We first recall
the definition of property $N_p$ of Green-Lazarsfeld. Let $X$ be a
smooth projective variety and $L \in \mbox{Pic}X$ a very ample
line bundle. Consider the embedding
\begin{equation*}
X \hookrightarrow \P H^0 (X,L) = \P
\end{equation*}
defined by the complete linear system of $L$. For the homogeneous
coordinate ring $S=\oplus_{\ell \geq 0} Sym^{\ell} H^0 (X,L)$ of
$\P$ and the finitely generated graded $S$-module $E=\oplus_{\ell
\in \Z} H^0 (X,L^{\ell})$, let
\begin{eqnarray*}
\cdots \rightarrow \oplus_j S^{k_{i,j}}(-i-j) \rightarrow  \cdots
\rightarrow \oplus_j S^{k_{1,j}}(-1-j) \rightarrow \oplus_j
S^{k_{0,j}}(-j) \rightarrow E \rightarrow 0
\end{eqnarray*}
be a minimal graded free resolution. Then $L$ is said to satisfy
property $N_p$ if $E$ admits a minimal free resolution of the form
\begin{eqnarray*}
\cdots \rightarrow S^{m_p}(-p-1) \rightarrow \cdots \rightarrow
S^{m_2}(-3) \rightarrow  S^{m_1}(-2) \rightarrow S \rightarrow E
\rightarrow 0.
\end{eqnarray*}
\noindent Therefore property $N_0$ holds if and only if $X
\hookrightarrow \P H^0 (X,L)$ is a projectively normal embedding,
property $N_1$ holds if and only if Property $N_0$ is satisfied
and the homogeneous ideal is generated by quadrics, and property
$N_p$ holds for $p \geq 2$ if and only if it has property $N_1$
and the $k^{th}$ syzygies among the quadrics are
generated by linear syzygies for all $1 \leq k \leq p-1$.\\

The aim of this article is to study higher syzygies of ruled
surfaces over irrational curves. More precisely we investigate
extremal curves on ruled surfaces. Recall that on a variety $X$
and a very ample line bundle $L \in \mbox{Pic}X$ which satisfies
property $N_p$ but not property $N_{p+1}$, a curve $C \subset X$
is said to be \textit{extremal} with respect to $L$ and property
$N_p$ if $(C,L_C)$ satisfies property $N_p$ but not property
$N_{p+1}$. See Remark 1.5 in \cite{GP4}. We will follow the
notation and terminology of R. Hartshorne's book \cite{H}, V $\S
2$. Let $C$ be a smooth projective curve of genus $g \geq 1$ and
let $\mathcal{E}$ be a vector bundle of rank $2$ on $C$ which is
normalized, i.e., $H^0 (C,\mathcal{E}) \neq 0$ while $H^0
(C,\mathcal{E} \otimes \mathcal{O}_C (D))=0$ for every divisor $D$
of negative degree. We set
\begin{equation*}
\mathfrak{e}=\wedge^2 \mathcal{E}~~~~\mbox{and}~~~~e = -
\mbox{deg}(\mathfrak{e}).
\end{equation*}
Let $X = \P_C (\mathcal{E})$ be the associated ruled surface with
projection morphism $\pi : X \rightarrow C$. We fix a minimal
section $C_0$ such that $\mathcal{O}_X (C_0)=\mathcal{O}_{\P_C
(\mathcal{E})} (1)$. For $\mathfrak{b} \in \mbox{Pic}C$,
$\mathfrak{b}f$ denote the pullback of $\mathfrak{b}$ by $\pi$.
Thus any element of $\mbox{Pic}X$ can be written
$aC_0+\mathfrak{b}f$ with $a\in \Z$ and $\mathfrak{b} \in
\mbox{Pic}C$ and any element of $\mbox{Num}X$ can be written $aC_0
+bf$ with $a,b \in \Z$.

When $X$ is an elliptic ruled surface with $e \geq 0$, the author
proved in \cite{ES} that a line bundle $L \in \mbox{Pic}X$ in the
numerical class $aC_0 + bf$ satisfies property $N_p$ if and only
if $a\geq1$ and $b-ae \geq 3+p$. Therefore when $L$ is very ample,
$L$ satisfies property $N_p$ if and only if $L_C$ satisfies
property $N_p$ where $L_C$ is the restriction of $L$ to a minimal
section $C_0$. In other words, for every very ample line bundle $L
\in \mbox{Pic}X$, $C_0$ is the extremal curve with respect to $L$
and property $N_p$. Also property $N_p$ is characterized in terms
of the intersection number of $L$ with a minimal section and a
fiber. Now we turn to the case when $X$ is an elliptic ruled
surface with $e=-1$. In \cite{GP2}, F. J. Gallego and B. P.
Purnaprajna conjectured the following:\\

\noindent {\bf Conjecture.} Let $X$ be an elliptic ruled surface
with $e=-1$ and $L \in \mbox{Pic}X$ a line bundle in the numerical
class $aC_0 + bf$. Then $L$ satisfies property $N_p$ if and only
if $a \geq 1$, $a+b \geq 3+p$, and $a+2b \geq 3+p$. \\

\noindent Note that there exists a smooth elliptic curve $E
\subset X$ such that $E \equiv 2C_0 -f$. See Proposition 3.2 in
\cite{GP1}. Also since
\begin{equation*}
\mbox{deg} (L|_{C_0}) = L \cdot C_0 = a+b~~\mbox{and}~~~\mbox{deg}
(L|_{E}) = L \cdot E = a+2b,
\end{equation*}
$(C,L|_{C_0})$ satisfies property $N_p$ if and only if $a+b \geq
3+p$, and $(E,L|_E)$ satisfies property $N_p$ if and only if $a+2b
\geq 3+p$. Therefore this conjecture suggests that when $L$ is
very ample, $(X,L)$ satisfies property $N_p$ if and only if
$(C,L|_{C_0})$ and $(E,L|_{E})$ satisfy property $N_p$, i.e.,
$C_0$ and $E$ are the extremal curves with respective to $L$ and
property $N_p$. This conjecture has been solved for $p=0$ by Yuko
Homma\cite{Homma2}, for $p=1$ by F. J. Gallego and B. P.
Purnaprajna\cite{GP1}, and for $a=1$ by the author\cite{ES}. And
our first result is the following:

\begin{theorem}\label{thm:elliptic}
Let $X$ be an elliptic ruled surface with $e=-1$ and $L \in
\mbox{Pic}X$ a line bundle in the numerical class $aC_0 + bf$. If
$L$ satisfies property $N_p$, then $a \geq 1$, $a+b \geq 3+p$ and
$a+2b \geq 3+p$.
\end{theorem}

\noindent Therefore the relevant part of the above
Gallego-Purnaprajna' conjecture is proved. For the proof we use D.
Eisenbud, M. Green, K. Hulek and S. Popescu's recent
work\cite{EGHP}. They investigate the relation between linear
syzygies in the minimal free resolution of a projective scheme $X
\subset \P^r$ and that of the linear sections of $X$. It is very
interesting that the failure of property $N_p$ for some $p$ is
closely related to the existence of a $(p+2)$-secant $p$-plane.
See Theorem \ref{thm:EGHP}. \\

\noindent {\bf Remark 1.1.} For $L \in \mbox{Pic}X$ in the
numerical class of $2C_0 + 3f$, property $N_2$ holds. See $\S 7$
in \cite{GP1}. This makes Gallego-Purnaprajna' conjecture
affirmative. Also note that although the conjecture is still open,
it is proved at Theorem 1.4.$(2)$ in \cite{ES} that for $L \in
\mbox{Pic}X$ in the
numerical class of $aC_0 + bf$, property $N_p$ holds if $a \geq 1$ and $a+2b \geq 5+2p$. \\

Now assume that $C$ is a curve of genus $g \geq 2$ and let $X$ be
a ruled surface over $C$ with $e \geq 0$. As in the case of
elliptic ruled surfaces with $e \geq 0$, one may expect that $C_0$
is the extremal curve with respect to $L$ and property $N_p$ for
every very ample line bundle $L \in \mbox{Pic}X$. And in this
article we prove that this is true for some cases. We first prove
the following:

\begin{theorem}\label{thm:main}
Let $X$ be a ruled surface over a curve $C$ of genus $g \geq 2$
with $e \geq 0$ and $L \in \mbox{Pic}X$ a line bundle in the
numerical class of $aC_0 +bf$ such that $a \geq 1$.\\
$(1)$ When $0 \leq e \leq g-3$, property $N_p$ holds for $L$ if
$b-ae \geq 3g-1-e+p$.\\
$(2)$ When $e \geq g-2$, property $N_p$ holds for $L$ if $b-ae
\geq 2g+1+p$.
\end{theorem}

\noindent This refines Theorem 1.8 in \cite{ES}. Also we prove the
following interesting relation between minimal free resolution of
a ruled surface and that of a minimal section.

\begin{theorem}\label{thm:minimalsection}
Let $X$ be a ruled surface over a curve $C$ of genus $g \geq 2$
with $e \geq 0$ and $L \in \mbox{Pic}X$ a line bundle in the
numerical class of $aC_0 +bf$ such that $a \geq 1$. Assume that
$b-ae \geq \mbox{max} \{2g+1, \frac{3}{2} g + p \}$ for some
nonnegative integer $p$. If $(X,L)$ satisfies property $N_p$, then
$(C,L_C)$ satisfies property $N_p$ where $L_C$ is the restriction
of $L$ to $C_0$.
\end{theorem}

\noindent In other words, $(X,L)$ fails to satisfy property $N_p$
if $(C,L_C)$ fails to satisfy property $N_p$. Therefore the
failure of the linearity of the minimal free resolution of $(X,L)$
is determined by that of $(C,L_C)$. For the proof, we use Theorem
(3.b.1) in \cite{Green} where M. Green studied the behavior of
Koszul cohomology under restriction to a divisor. Roughly
speaking, the relation between higher syzygies of $X$ and that of
$C_0$ is encoded in the short exact sequence
\begin{equation*}
0 \rightarrow \mathcal{O}_X (-C_0) \rightarrow \mathcal{O}_X
\rightarrow \mathcal{O}_C \rightarrow 0.
\end{equation*}
More precisely, if $e \geq 0$ and if $L \in \mbox{Pic}X$ is in the
numerical class $aC_0 +bf$ such that $a \geq 1$ and $b-ae \geq
2g+1$, then $H^1 (X,L^j \otimes \mathcal{O}_X (-C_0))=0$ for all
$j \geq 0$. Thus we have the short exact sequence
\begin{equation*}
0 \rightarrow R_1 \rightarrow R_2 \rightarrow R_3 \rightarrow 0
\end{equation*}
of graded $S$-modules with maps preserving the gradings where
\begin{equation*}
R_1 =\oplus_{\ell \in \Z} H^0 (X,\mathcal{O}_X (-C_0) \otimes
L^{\ell}),~R_2 =\oplus_{\ell \in \Z} H^0
(X,L^{\ell})~\mbox{and}~R_3 =\oplus_{\ell \in \Z} H^0 (C,L_C
^{\ell})
\end{equation*}
and $L_C$ is the restriction of $L$ to $C_0$. This gives a long
exact sequence for Koszul cohomology groups which enables us to
compare the minimal free resolution of $R_2$ and that of $R_3$.
Combining Theorem \ref{thm:main} and Theorem
\ref{thm:minimalsection}, we obtain
some corollaries. First recall the following Green-Lazarsfeld's result\cite{GL2}:\\

\begin{enumerate}
\item[$(i)$] Let $C$ be a curve of genus $g \geq 2$. For $L \in
\mbox{Pic}C$ with $\mbox{deg}(L) \geq 3g-2$, $L$ satisfies
property $N_p$ if and only if $\mbox{deg}(L) \geq 2g+1+p$.

\item[$(ii)$] If $C$ is a hyperelliptic curve of genus $g \geq 2$
and if $L \in \mbox{Pic}C$ is a very ample line bundle, then $L$
satisfies property $N_p$ if and only if $\mbox{deg}(L) \geq
2g+1+p$.\\
\end{enumerate}
We generalize this as follows:

\begin{corollary}\label{cor:arbitrary}
Let $X$ be a ruled surface over a curve $C$ of genus $g \geq 3$
with $e \geq 0$ and $L \in \mbox{Pic}X$ a line bundle in the
numerical class of $aC_0 +bf$ such that $a \geq 1$. Assume that
$b-ae \geq 3g-2$.\\
$(1)$ If $L$ satisfies property $N_p$, then $b-ae \geq 2g+1+p$.\\
$(2)$ When $e \geq g-2$, $L$ satisfies property  $N_p$ if and only
if $b-ae \geq 2g+1+p$.
\end{corollary}

\begin{theorem}\label{thm:hyperellipticnormalgeneration}
Let $C$ be a hyperelliptic curve of genus $g \geq 2$ and $X$ a
ruled surface over $C$ with $e \geq 0$. Let $L \in \mbox{Pic}X$be
a very ample line bundle  in the numerical class of $aC_0 +bf$.\\
$(1)$ $L$ is normally generated if and only if $b-ae \geq 2g+1$.\\
$(2)$ $L$ is normally presented if and only if $b-ae \geq 2g+2$.\\
$(3)$ If $L$ satisfies property $N_p$, then $b-ae \geq 2g+1+p$.\\
$(4)$ Assume that $e \geq g-2$. Then $L$ satisfies property $N_p$
if and only if $b-ae \geq 2g+1+p$.
\end{theorem}

\noindent {\bf Remark 1.2.} When $C$ is a curve of genus $g=2$,
let $X$ be a ruled surface over $C$ with $e \geq 0$. Then by
Corollary \ref{cor:genustwo}, $L \in \mbox{Pic}X$ in the numerical
class of $aC_0 +bf$ is very ample if and only if $a \geq 1$ and
$b-ae \geq 5$. Also Theorem
\ref{thm:hyperellipticnormalgeneration} says that $L$ satisfies
property $N_p$ if and only if $L_C$ satisfies property $N_p$.
Therefore as in the case of elliptic ruled surfaces with $e \geq
0$, a minimal section $C_0$ is the
extremal curve with respect to $L$ and property $N_p$. \qed \\

The organization of this paper is as follows. In $\S 2$, we review
some elementary facts to study higher syzygies of ruled surfaces.
$\S 3$ is devoted to prove the relevant part of
Gallego-Purnaprajna' conjecture. In $\S 4$, we develop the
technique to study higher syzygies of ruled surfaces with $e \geq
0$ and prove some new results. Finally $\S 5$ is devoted to apply
our results to ruled surfaces over a hyperelliptic curve with $e
\geq 0$.\\

\section{Preliminary}
\subsection{Notations and Conventions}
Throughout this paper the following is assumed.\\

\noindent $(1)$ All varieties are defined over the complex number
field $\C$.\\

\noindent $(2)$ For a finite dimensional $\C$-vector space $V$,
$\P(V)$ is the projective space of one-

dimensional quotients of $V$.\\

\noindent $(3)$ When a variety $X$ is embedded in a projective
space, we always assume that

it is non-degenerate, i.e. it does not lie
in any hyperplane.\\

\noindent $(4)$ When a projective variety $X$ is embedded in a
projective space $\P^r$ by a very ample

line bundle $L \in \mbox{Pic} X $, we may write $\mathcal{O}_X
(1)$ instead of $L$ as long as there is no

confusion.\\

\subsection{Minimal slope} Let $C$ be a smooth projective curve of genus $g \geq 1$. For a vector
bundle $\mathcal{F}$ on $C$, the slope $\mu(\mathcal{F})$ is
defined by $\deg(\mathcal{F}) / \mbox{rank}(\mathcal{F})$ and the
minimal slope $\mu^- (\mathcal{F})$ is defined as follows:
\begin{equation*}
\mu^- (\mathcal{F})=\mbox{min}\{ \mu(Q) | \mathcal{F}\rightarrow Q
\rightarrow 0\}
\end{equation*}
It is well-known that $\mu^- (\mathcal{F})$ measures the
positivity of $\mathcal{F}$.

\begin{lemma}\label{lem:folklore}
$(1)$ If $\mu^- (\mathcal{F}) > 2g-2$, then $h^1 (C,
\mathcal{F})=0$.\\
$(2)$ If $\mu^- (\mathcal{F})  > 2g-1$, then $\mathcal{F}$ is
globally generated.\\
$(3)$ If $\mu^- (\mathcal{F})  > 2g$, then
$\mathcal{O}_{\P(\mathcal{F})} (1)$ is very ample.
\end{lemma}

\begin{proof}
See $\S 1$ in \cite{Butler}.
\end{proof}

\begin{lemma}\label{lem:veryampleness}
For a ruled surface $X$ over $C$ with $e$, let $L \in \mbox{Pic}X$
be a line bundle in the numerical class of $aC_0 +bf$ with $a \geq 1$
and put $\mathcal{F} = \pi_* L$.\\
$(1)$ If $\mathcal{O}_{\P(\mathcal{F})} (1)$ is very ample, then
$L$ is very ample.\\
$(2)$ If $\mu^- (\mathcal{F}) > 2g$, then $L$ is very ample.
\end{lemma}

\begin{proof}
$(1)$ Note that $X \subset \P (\mathcal{F})$ is given by a
fiberwise $a$-uple map and $\mathcal{O}_{\P (\mathcal{F})}
(1)|_{X} = L$. Since $\mathcal{O}_{\P (\mathcal{F})} (1)$ is very
ample, $L$ is very ample.\\
$(2)$ By Lemma \ref{lem:folklore}.$(3)$,
$\mathcal{O}_{\P(\mathcal{F})} (1)$ is very ample. Therefore $L$
is very ample by $(1)$.
\end{proof}

\noindent For elliptic ruled surfaces, very ampleness of line
bundles is determined by its numerical type. And the following
shows that this is true for ruled surfaces over a curve $C$ of
genus $2$ with $e \geq 0$.

\begin{corollary}\label{cor:genustwo}
For a curve $C$ of genus $2$, let $X$ be a ruled surface over $C$
with $e \geq 0$. Then $L \in \mbox{Pic}X$ in the numerical class
of $aC_0 +bf$ is very ample if and only if $a \geq 1$ and $b-ae
\geq 5$.
\end{corollary}

\begin{proof}
If $L$ is very ample, then the restrictions of $L$ to a fiber and
to a minimal section $C_0$ are both very ample. Thus $a \geq 1$
and $\mbox{deg}(L|_{C_0})=b-ae \geq 5$. Conversely, if $a \geq 1$
and $\mu^- (\pi_* L)=b-ae \geq 5$, then $L$ is very ample by Lemma
\ref{lem:veryampleness}.$(2)$.
\end{proof}

\subsection{Cohomological interpretation of property $N_p$} We
review some cohomological criteria for property $N_p$. Let $X$ be
a smooth projective variety of dimension $n \geq 1$ and let $L \in
\mbox{Pic}X$ be a very ample line bundle. Consider the natural
short exact sequence
\begin{equation*}
0 \rightarrow \mathcal{M}_L \rightarrow H^0 (X,L) \otimes
\mathcal{O}_X \rightarrow L \rightarrow 0.
\end{equation*}

\noindent It is well-known that $L$ satisfies property $N_p$ if
and only if $H^1 (X,\wedge^i \mathcal{M}_L \otimes L^j)=0$ for $1
\leq i \leq p+1$ and all $j \geq 1$. And in the situation that
will concern us, one can get away with checking even a little
less:

\begin{lemma}\label{lem:criterion}
Suppose that the ideal sheaf $~\mathcal{I}_{X/\P}~$ of $X
\hookrightarrow \P = \P H^0 (X,L)$ is $3$-regular in the sense of
Castelnuovo-Mumford, i.e., that $H^i (\P,\mathcal{I}_{X/\P}
(3-i))=0$ for all $i \geq 1$. Then for $p \leq
\mbox{codim}(X,\P)$, property $N_p$ holds for $L$ if and only if
\begin{equation*}
H^1 (X,\wedge^{p+1} \mathcal{M}_L \otimes L)=0.
\end{equation*}
\end{lemma}

\begin{proof}
See $\S 1$ in \cite{GL2}.
\end{proof}

\noindent {\bf Remark 2.1.} Let $X$ be a smooth projective surface
with geometric genus $0$, i.e., $H^2 (X,\mathcal{O}_X)=0$. Let $L
\in \mbox{Pic}X$ be a normally generated very ample line bundle
such that $H^1 (X,L)=0$.
Then  it is easy to check that $X \hookrightarrow \P H^0 (X,L)=\P$, $\mathcal{I}_{X/\P}$ is $3$-regular. \qed \\

\subsection{Higher syzygies of degenerate varieties}
Let $\Lambda \cong \P W \subset \P V$ be a linear subspace such
that $\mbox{codim}(W,V)=c$. It is easily checked that
\begin{equation*}
\Omega_{\P V} (1)|_{\P W} \cong \Omega_{\P W} (1) \oplus
\mathcal{O}_{\P W} ^c.
\end{equation*}
Now let $X \subset \P V$ be a smooth projective variety which is
indeed contained in $\Lambda$. Let the corresponding very ample
line bundle on $X$ be $L \in \mbox{Pic}X$. Consider the natural
short exact sequences
\begin{eqnarray*}
& 0 \rightarrow \mathcal{M}_V  \rightarrow V \otimes \mathcal{O}_X \rightarrow L \rightarrow  0&~~~~\mbox{and} \\
& 0 \rightarrow \mathcal{M}_W \rightarrow W \otimes \mathcal{O}_X
\rightarrow L \rightarrow  0.
\end{eqnarray*}
The above observation shows that $\mathcal{M}_V \cong
\mathcal{M}_W \oplus \mathcal{O}_{X} ^c$.

\begin{lemma}\label{lem:degenrate}
Under the situation just stated, assume that $H^1 (X,L^{\ell}) =0$
for all $\ell \geq 1$ and $W = H^0 (X,L)$. Then $(X,L)$ satisfies
property $N_p$ if and only if
\begin{equation*}
H^1 (X,\wedge^i \mathcal{M}_V \otimes L^{\ell}) = 0~~~~\mbox{for
$1 \leq i \leq p+1$ and $\ell \geq 1$.}
\end{equation*}
\end{lemma}

\begin{proof}
See Lemma 2.6 in \cite{ES}.
\end{proof}

\section{Elliptic ruled surfaces with $e=-1$}
\noindent This section is devoted to prove Theorem
\ref{thm:elliptic}. We first remark the following recent result:

\begin{theorem}[D. Eisenbud, M. Green, K. Hulek and S. Popescu, Theorem 1.1 in
\cite{EGHP}]\label{thm:EGHP} For a projective variety $X$ and a
very ample line bundle $L \in \mbox{Pic}X$, if $X \subset \P H^0
(X,L)$ admits a $(p+2)$-secant $p$-plane, i.e., there exists a
linear subspace $\Lambda \subset \P H^0 (X,L)$ of dimension $\leq
p$ such that $X \cap \Lambda$ is finite and $\mbox{length}(X \cap
\Lambda) \geq p+2$, then $(X,L)$ fails to satisfy Property $N_p$.
\end{theorem}

\noindent For elliptic ruled surfaces with $e=-1$, Theorem
\ref{thm:EGHP} enables us to relate the failure of property $N_p$
to the minimal free resolution of the minimal section $C_0$ and
that of the elliptic curve $E \equiv 2C_0 -f$. See the following proof:\\

\noindent {\bf Proof of Theorem \ref{thm:elliptic}.} Assume that
$L$ satisfies property $N_p$. When $p=0$, Yuko Homma\cite{Homma2}
proved that $L$ satisfies property $N_0$ if and only if $a \geq
1$, $a+b \geq 3$ and $a+2b \geq 3$. When $p=1$, F. J. Gallego and
B. P. Purnaprajna\cite{GP1} proved that $L$ satisfies property
$N_1$ if and only if $a \geq 1$, $a+b \geq 4$ and $a+2b \geq 4$.
Thus we assume that $p \geq 2$ and hence $a \geq 1$, $a+b \geq 4$
and $a+2b \geq 4$.

We first prove that $a+b \geq 3+p$. Assume that $a+b = 2+p(p \geq
2)$. Note that $\mbox{deg}(L|_{C_0})=a+b$. Also $H^0 (X,L)
\rightarrow H^0 (C,L|_{C_0})$ is surjective since $H^1 (X,L
\otimes \mathcal{O}_X (-C_0))=0$. See Proposition 4.3 in
\cite{GP2}. Thus $C_0$ is embedded in a linear subspace $\Lambda
\cong \P^{p+1} \subset \P H^0 (X,L)$ by $X \hookrightarrow \P H^0
(X,L)$. Let $\Lambda' \cong \P^p \subset \Lambda$ be a hyperplane
such that $\Lambda' \cap X$ is a finite scheme. Then
\begin{equation*}
\mbox{length}(X \cap \Lambda') \geq \mbox{length}(C_0 \cap
\Lambda') = p+2.
\end{equation*}
Thus Theorem \ref{thm:EGHP} implies that $X \hookrightarrow \P H^0
(X,L)$ fails to satisfy property $N_p$ which contradicts to our
assumption. Therefore $a+b \geq 3+p$.

By the same way we prove that $a+2b \geq 3+p$. Assume that $a+2b =
2+p(p \geq 2)$. Note that there exists a smooth elliptic curve $E
\subset X$ such that $E \equiv 2C_0 -f$ and
$\mbox{deg}(L|_E)=a+2b$. Also $H^0 (X,L) \rightarrow H^0 (E,L|_E)$
is surjective since $H^1 (X,L \otimes \mathcal{O}_X (-E))=0$ by
Proposition 4.3 in \cite{GP2}. Thus $E$ is embedded in a linear
subspace $\Lambda \cong \P^{p+1} \subset \P H^0 (X,L)$ by $X
\hookrightarrow \P H^0 (X,L)$. Let $\Lambda' \cong \P^p \subset
\Lambda$ be a hyperplane such that $\Lambda' \cap X$ is a finite
scheme. Then
\begin{equation*}
\mbox{length}(X \cap \Lambda') \geq \mbox{length}(E \cap \Lambda')
= p+2.
\end{equation*}
Thus $(X,L)$ cannot satisfy property $N_p$ by Theorem
\ref{thm:EGHP}.  Therefore $a+2b \geq 3+p$. \qed \\

\noindent {\bf Remark 3.1.} For a smooth projective variety $X$ of
dimension $n$ and an ample line bundle $A \in \mbox{Pic}X$, it is
conjectured by Shigeru Mukai that $K_X + (n+2+p)A$ satisfies
property $N_p$. Although this is still open, Theorem
\ref{thm:elliptic} shows that the condition of the conjecture is
optimal when $n=2$. Indeed let $X$ be an elliptic ruled surface
with $e =-1$ and let $A \in \mbox{Pic}X$ be a line bundle in the
numerical class $aC_0 + bf$ such that $a \geq 1$ and $a+b = 1$.
Then
\begin{equation*}
K_X + (4+p)A \equiv \{ a(4+p)-2 \}C_0 + \{ b(4+p)+1 \}f
\end{equation*}
and hence $\{ a(4+p)-2 \} + \{ b(4+p)+1 \} = (4+p)(a+b)-1=3+p$.
Thus Theorem \ref{thm:elliptic} says that $K_X + (4+p)A$ fails to
satisfy property $N_{p+1}$. For a more general result, see $\S 5$
in \cite{ES2}. Indeed for all $n \geq 2$ and $p \geq 0$  it is
proved that there exists an $n$-dimensional ruled variety $X$ over
a curve and an ample line bundle $A \in \mbox{Pic}X$ such that for
any $p \geq 0$, $(X,K_X + (n+2+p)A)$ fails to satisfy Property
$N_{p+1}$. \qed \\

\section{Ruled surfaces with $e \geq 0$}
\noindent Let $C$ be a smooth curve of genus $g \geq 2$ and $X$ a
ruled surface over $C$ with $e \geq 0$. For a line bundle $L \in
\mbox{Pic}X$ in the numerical class of $aC_0 +bf$, assume that $a
\geq 1$ and $b-ae \geq 2g+1$. Then $L$ is very ample and $H^1
(X,L^j)=0$ for all $j \geq 1$. Also the vector bundle
$\mathcal{F}=\pi_* L$ over $C$ is globally generated and $\mu^-
(\mathcal{F})=b-ae$. Therefore we have the following two short
exact sequences:
\begin{eqnarray*}
0 \rightarrow  \mathcal{M}_{\mathcal{F}} \rightarrow  H^0
(C,\mathcal{F}) \otimes \mathcal{O}_C
\rightarrow  \mathcal{F} \rightarrow 0 &~~~~\mbox{and} \\
0 \rightarrow  \mathcal{M}_L \rightarrow H^0 (X,L) \otimes
\mathcal{O}_X \rightarrow L \rightarrow 0.\\
\end{eqnarray*}

\noindent We first prove the following simple criterion for
property $N_p$:

\begin{lemma}\label{lem:test}
Under the situation just stated, assume that $b-ae=2g+1+q$ for
some $q \geq 0$. Then for $p \leq  \frac{g}{2} +1+q$, $(X,L)$
satisfies property $N_p$ if and only if
\begin{equation*}
H^1 (X,\wedge^{p+1} M_L \otimes  L)=0.
\end{equation*}
\end{lemma}

\begin{proof}
Since $b-ae \geq 2g+1$, $L$ is normally generated by Theorem 5.1A
in \cite{Butler}. Also $H^1 (X,L)=H^2 (X,\mathcal{O}_X)=0$.
Therefore $X \hookrightarrow \P H^0 (X,L)$ is $3$-regular in the
sense of Castelnuovo-Mumford. Thus for $p \leq h^0 (X,L)-3$, $L$
satisfies property $N_p$ if and only if $H^1 (X,\wedge^{p+1} M_L
\otimes  L)=0$ by Lemma \ref{lem:criterion} or by Lemma 1.10 in
\cite{GL2}. Also by Riemann-Roch,
\begin{equation*}
h^0 (X,L)-3= a(g+2+q)+g +q-1+ \frac{ae(a+1)}{2} \geq \frac{g}{2}
+1+q,
\end{equation*}
which completes the proof.
\end{proof}

\begin{lemma}\label{lem:generalvanishing}
Under the situation just stated, let $p \geq 0$ be a given
integer.\\
$(1)$ $H^1 (X,\wedge^{p+1} M_L \otimes L \otimes \mathcal{O}_X
(-C_0))=0$ if $$ b-ae \geq \begin{cases} \mbox{max}\{2g+1,3g-1-e+p\} & \mbox{when $0\leq e \leq 2g-2$,} \\
\mbox{max}\{2g+1,g+1+p\} & \mbox{when $2g-1\leq e$, and} \\
\mbox{max}\{2g+1,g+p\} & \mbox{when $2g \leq e$.}
\end{cases}$$
$(2)$ $H^1 (X,\wedge^p M_L \otimes L^2 \otimes \mathcal{O}_X
(-C_0))=0$ if $b-ae \geq \mbox{max}\{2g+1,\frac{3}{2}g+p \}$.
\end{lemma}

\begin{proof}
Consider the exact sequence
\begin{eqnarray*}
0 \rightarrow  \mathcal{K}_L \rightarrow \pi^* \mathcal{F}
\rightarrow L \rightarrow 0
\end{eqnarray*}
where $\mathcal{K}_L$ is a vector bundle of rank $a$ on $X$ which
is $1$ $\pi$-regular. Using Snake Lemma, we have the following
commutative diagram:
\begin{equation*}
\begin{CD}
&  &&  && 0 &\\
&  &&  && \downarrow &\\
& 0 && && \mathcal{K}_L &\\
& \downarrow &&  && \downarrow &\\
0 \rightarrow &\pi^* \mathcal{M}_\mathcal{F}& \rightarrow & H^0 (C,\mathcal{F}) \otimes \mathcal{O}_X & \rightarrow  & \pi^* \mathcal{F} & \rightarrow 0 \\
& \downarrow && \parallel \wr && \downarrow &\\
0 \rightarrow & \mathcal{M}_L & \rightarrow & H^0 (X,L) \otimes \mathcal{O}_X & \rightarrow & L & \rightarrow 0 \\
& \downarrow &&  && \downarrow &\\
& \mathcal{K}_L && && 0  &\\
& \downarrow  &&  && &\\
& 0 &&  && &.\\\\
\end{CD}
\end{equation*}
$(1)$ From the sequence
\begin{equation*}
0 \rightarrow \pi^* M_\mathcal{F} \rightarrow M_L \rightarrow
\mathcal{K}_{L} \rightarrow 0,
\end{equation*}
$H^1 (X,\wedge^{p+1} M_L \otimes L \otimes \mathcal{O}_X
(-C_0))=0$ if
\begin{equation*}
H^1 (X,\wedge^{p+1-i} \pi^* M_{\mathcal{F}} \otimes \wedge^{i}
\mathcal{K}_{L}  \otimes L \otimes \mathcal{O}_X (-C_0))=0
\end{equation*}
for every $0 \leq i \leq \mbox{min} \{a,m\}$.\\

\noindent \underline{\textit{Case 1.}}\quad First we concentrate
on the case $i =0$. Then
\begin{equation*}
H^1 (X,\wedge^{p+1} \pi^*M_{\mathcal{F}} \otimes L \otimes
\mathcal{O}_X (-C_0))\cong H^1 (C,\wedge^{p+1} M_{\mathcal{F}}
\otimes \pi_* \{L \otimes \mathcal{O}_X (-C_0)\})
\end{equation*}
and hence it suffices to show that
\begin{equation*}
\mu^- (\wedge^{p+1} M_{\mathcal{F}} \otimes \pi_* \{L \otimes
\mathcal{O}_X (-C_0)\}) > 2g-2.
\end{equation*}
Indeed we have
\begin{eqnarray*}
\mu^- (\wedge^{p+1} M_{\mathcal{F}} \otimes \pi_* \{L \otimes
\mathcal{O}_X (-C_0)\}) & \geq & (p+1) \mu^- (M_{\mathcal{F}})
+  \mu^- ( \mathcal{F})+e\\
& \geq & - (p+1) \frac{\mu^- (\mathcal{F})}{\mu^- (\mathcal{F})-g}
+ \mu^- (\mathcal{F}) +e
\end{eqnarray*}
since $\mu^- (\mathcal{E}) =e$ when $e \geq 0$. Therefore it
suffices to show that
\begin{equation*}
- (p+1) \frac{\mu^- (\mathcal{F})}{\mu^- (\mathcal{F})-g} + \mu^-
(\mathcal{F}) +e >2g-2
\end{equation*}
or equivalently
\begin{equation*}
\mu^- (\mathcal{F})^2 + (e+1-3g-p)\mu^- (\mathcal{F}) + g(2g-2-e) >0.
\end{equation*}
It is a tedious calculation that the second inequality holds under our assumption.\\

\noindent \underline{\textit{Case 2.}}\quad Now we consider the
case $1 \leq i \leq a-1$.  Since $\wedge^{i} \mathcal{K}_{L}$ is a
direct summand of the tensor product $T^{i} \mathcal{K}_{L}$, it
suffices to show that
\begin{equation*}
H^1 (X,\wedge^{p+1-i} \pi^* M_{\mathcal{F}} \otimes T^{i}
\mathcal{K}_{L} \otimes L \otimes \mathcal{O}_X (-C_0))=0.
\end{equation*}
Note that since $i \leq a-1$, $T^{i} \mathcal{K}_{L} \otimes L
\otimes \mathcal{O}_X (-C_0)$ is $0$ $\pi$-regular. Therefore
\begin{equation*}
H^1 (X,\wedge^{p+1-i} \pi^* M_{\mathcal{F}} \otimes T^{i}
\mathcal{K}_{L} \otimes L \otimes \mathcal{O}_X (-C_0) )=H^1
(C,\wedge^{p+1-i} M_{\mathcal{F}} \otimes \pi_* \{T^{i}
\mathcal{K}_{L} \otimes L \otimes \mathcal{O}_X (-C_0)\}).
\end{equation*}
Note that
\begin{eqnarray*}
\mu^- (\wedge^{p+1-i} M_{\mathcal{F}} \otimes \pi_* \{T^{i}
\mathcal{K}_{L} \otimes L \otimes \mathcal{O}_X (-C_0)\}) & \geq &
- (p+1-i) \frac{\mu^- (\mathcal{F})}{\mu^- (\mathcal{F})-g} +
(i+1) \mu^- (\mathcal{F}) + e \\
& \geq & - (p+1) \frac{ \mu^- (\mathcal{F})}{\mu^-
(\mathcal{F})-g} +  \mu^- (\mathcal{F}) + e \\
& > & 2g-2
\end{eqnarray*}
by the claim in the proof of Lemma 3.1.$(1)$ in \cite{ES}.\\

\noindent \underline{\textit{Case 3.}}\quad Assume that $i=a$.
Note that $\wedge^a \mathcal{K}_L = \wedge^{a+1} \pi^* \mathcal{F}
\otimes L^{-1}$. Using this equality, we have
\begin{eqnarray*}
\wedge^{p+1-a} \pi^* M_{\mathcal{F}} \otimes \wedge^{a}
\mathcal{K}_{L} \otimes L \otimes \mathcal{O}_X (-C_0)
=\wedge^{p+1-a} \pi^* M_{\mathcal{F}} \otimes \wedge^{a+1} \pi^*
\mathcal{F} \otimes \mathcal{O}_X (-C_0)
\end{eqnarray*}
and hence this is $1$ $\pi$-regular. Therefore
\begin{eqnarray*}
& & H^1 (X,\wedge^{p+1-a} \pi^* M_{\mathcal{F}} \otimes \wedge^{a}
\mathcal{K}_{L} \otimes \mathcal{O}_X (-C_0) )   \\
& = & H^1 (C,\wedge^{p+1-a}  M_{\mathcal{F}} \otimes \wedge^{a+1}
\mathcal{F} \otimes \pi_* \mathcal{O}_X (-C_0))=0
\end{eqnarray*}
since $\pi_*  \mathcal{O}_X (-C_0) =0$. \\\\
$(2)$ We apply Lemma 3.1.$(1)$ in \cite{ES} to $L^2 \otimes
\mathcal{O}_X (-C_0) \equiv (2a-1)C_0 + 2bf$. Therefore we get the
desired vanishing if
\begin{equation*}
2 \mu^- (\mathcal{F})+e
> \frac{p\mu^- (\mathcal{F})}{\mu^- (\mathcal{F})-g}+2g-2.
\end{equation*}
It is easily checked that this inequality holds if $\mu^-
(\mathcal{F}) \geq \mbox{max}\{2g+1,\frac{3}{2}g+p \}$.
\end{proof}

\noindent Consider the short exact sequence
\begin{equation*}
0 \rightarrow \mathcal{O}_X (-C_0) \rightarrow \mathcal{O}_X
\rightarrow \mathcal{O}_C \rightarrow 0
\end{equation*}
and define graded $S$-modules
\begin{equation*}
R_1 =\oplus_{\ell \in \Z} H^0 (X,\mathcal{O}_X (-C_0) \otimes
L^{\ell}),~R_2 =\oplus_{\ell \in \Z} H^0
(X,L^{\ell})~\mbox{and}~R_3 =\oplus_{\ell \in \Z} H^0 (C,L_C
^{\ell})
\end{equation*}
where $L_C$ is the restriction of $L$ to $C_0$. Since $H^1
(X,\mathcal{O}_X (-C_0) \otimes L^{\ell})=0$ for all $\ell \geq
1$, we have the short exact sequence
\begin{equation*}
0 \rightarrow R_1 \rightarrow R_2 \rightarrow R_3 \rightarrow 0
\end{equation*}
of graded $S$-modules with maps preserving the gradings and hence
there is a long exact sequence
\begin{eqnarray*}
\cdots \rightarrow M_{1,q}(R_1,V) \rightarrow
M_{1,q}(R_2,V) \rightarrow M_{1,q}(R_3,V)\\
\rightarrow M_{0,q}(R_1,V) \rightarrow M_{0,q}(R_2,V) \rightarrow
M_{0,q}(R_3,V)\rightarrow 0.
\end{eqnarray*}
See Corollary (1.d.4) in \cite{Green}. We need the following part:
\begin{eqnarray*}
\cdots \rightarrow M_{p,p+2}(R_1,V) \rightarrow M_{p,p+2}(R_2,V)
\rightarrow M_{p,p+2}(R_3,V) \rightarrow M_{p-1,p+2}(R_1,V)
\rightarrow \cdots.
\end{eqnarray*}
Since $H^1 (X,\mathcal{O}_X (-C_0) \otimes L)=H^1 (X,\mathcal{O}_X
(-C_0) \otimes L^2)=0$,
\begin{eqnarray*}
& & M_{p,p+2}(R_1,V) \cong H^1 (X,\wedge^{p+1} \mathcal{M}_L
\otimes \mathcal{O}_X (-C_0) \otimes L)~\mbox{and}\\
& & M_{p-1,p+2}(R_1,V) \cong H^1 (X,\wedge^{p} \mathcal{M}_L
\otimes \mathcal{O}_X (-C_0) \otimes L^2).
\end{eqnarray*}
See Theorem 2.4 in \cite{ES}. Now we start to prove Theorem
\ref{thm:main} and Theorem \ref{thm:minimalsection}.\\

\noindent {\bf Proof of Theorem \ref{thm:main}.} First note that
if $\mbox{deg}(L_C)=b-ae \geq 2g+1+p$, then $(C,L_C)$ satisfies
property $N_p$. Therefore $M_{p,p+2}(R_3,V)\cong H^1
(C,\wedge^{p+1} \mathcal{M}_L \otimes L_C)=0$ by Lemma 2.6 in
\cite{ES}. Also $H^1 (X,\wedge^{p+1} \mathcal{M}_L \otimes
\mathcal{O}_X (-C_0) \otimes L)=0$ if
$$ b-ae \geq \begin{cases}
\mbox{max}\{2g+1,3g-1-e+p\} & \mbox{when $0\leq e \leq g-3$, and} \\
 2g+1 +p  & \mbox{when $e \geq g-2$.}
\end{cases}$$
by Lemma \ref{lem:generalvanishing}.$(1)$. Thus it is proved that
$M_{p,p+2}(R_2,V)\cong H^1 (X,\wedge^{p+1} \mathcal{M}_L \otimes
L) =0$. Therefore $(X,L)$ satisfies property $N_p$ by
Lemma \ref{lem:test}. \qed \\

\noindent {\bf Proof of Theorem \ref{thm:minimalsection}.} By
Lemma \ref{lem:generalvanishing}.$(2)$,
\begin{equation*}
M_{p-1,p+2}(R_1,V) \cong H^1 (X,\wedge^{p} \mathcal{M}_L \otimes
\mathcal{O}_X (-C_0) \otimes L^2)=0
\end{equation*}
if $b-ae \geq \mbox{max}\{2g+1,\frac{3}{2}g+p \}$. Also if $L$
satisfies property $N_{p}$, then $M_{p,p+2}(R_2,V)=0$. Thus
$M_{p,p+2}(R_3,V)=0$. Note that $H^1 (C,L_C)=0$ since
$\mbox{deg}(L_C)=b-ae \geq 2g+1$ and hence
\begin{equation*}
M_{p,p+2}(R_3,V) \cong H^1 (C,\wedge^{p+1} \mathcal{M}_L \otimes L
\otimes \mathcal{O}_C).
\end{equation*}
Therefore $(C,L_C)$ satisfies property $N_{p}$ by Lemma 2.6 in \cite{ES}. \qed \\

\noindent {\bf Proof of Corollary \ref{cor:arbitrary}.} $(1)$
Since $g \geq 3$, $b-ae \geq 3g-2 \geq 2g+1$. Assume that $b-ae =
2g+p$. If $(X,L)$ satisfies property $N_p$, then $(C,L_C)$
satisfies property $N_p$ by Theorem \ref{thm:minimalsection}. But
$H^0 (C,L_C-K_C) \neq 0$ since $\mbox{deg}(L_C)=b-ae=2g+p \geq
3g-2$. Thus M. Green and R. Lazarsfeld's result in \cite{GL2} says
that $(C,L_C)$ fails to satisfy property $N_p$ which is a
contradiction. Therefore if $L$ satisfies property $N_p$, then
$b-ae \geq 2g+1+p$.\\
$(2)$ This is proved immediately by $(1)$ and Theorem
\ref{thm:main}.$(2)$. \qed \\

\section{Ruled surfaces over a hyperelliptic curve with $e \geq 0$}
\noindent The main goal of this section is to prove Theorem
\ref{thm:hyperellipticnormalgeneration}. We first prove the
following criterion for normal generation:

\begin{proposition}\label{prop:normalgenerationsurface}
Let $S$ be a surface and let $C \subset S$ be a curve. For a
normally generated line bundle $L \in \mbox{Pic}S$, assume that
$H^1 (S,\mathcal{O}_S (-C) \otimes L^{j+1})=0$ for all $j \geq 1$.
Then $(C,L|_C)$ is normally generated.
\end{proposition}

\begin{proof}
For $j \geq 1$, consider the following commutative diagram:
\begin{equation*}
\begin{CD}
H^0 (X,L) \otimes H^0 (X,L^j) & \stackrel{\gamma}{\rightarrow} & H^0 (C,L_C) \otimes H^0 (C,L_C ^j)& \\
\downarrow  \alpha            &                                & \downarrow  \beta                  & \\
H^0 (X,L^{j+1})               & \stackrel{\delta}{\rightarrow} & H^0 (C,L_C ^{j+1})    &  \rightarrow H^1 (S,\mathcal{O}_S (-C) \otimes L^{j+1})=0        \\\\
\end{CD}
\end{equation*}
$\alpha$ is surjective since $L$ is normally generated. Also from
the short exact sequence
\begin{equation*}
0 \rightarrow \mathcal{O}_S (-C) \rightarrow \mathcal{O}_S
\rightarrow \mathcal{O}_C \rightarrow 0,
\end{equation*}
$\delta$ is surjective since $H^1 (S,\mathcal{O}_S (-C) \otimes
L^{j+1})=0$ for all $j \geq 1$. Therefore $\beta$ is also
surjective for all $j \geq 1$.
\end{proof}

\noindent {\bf Proof of Theorem
\ref{thm:hyperellipticnormalgeneration}.} It is well-known that
$L$ is normally generated if $b-ae \geq 2g+1$ and normally
presented if $b-ae \geq 2g+2$. See Theorem 5.1A and Theorem 5.1B
in \cite{Butler} or Corollary 4.6 in \cite{ES2}.\\
$(1)$ Assume that $L$ is normally generated. Note that for $L_C :=
L|_{C_0}$, $\mbox{deg}(L_C) =b-ae \geq g+3$ since $C$ is
hyperelliptic and $L_C$ is very ample. Thus one can easily check
that
\begin{equation*}
H^1 (X,L^{j+1} \otimes \mathcal{O}_X (-C_0))=0~~~\mbox{for
all}~~~~j \geq 1.
\end{equation*}
By Proposition \ref{prop:normalgenerationsurface}, $L_C$ is
normally generated and hence $b-ae \geq 2g+1$.\\
$(2)$ Assume that $L$ is normally presented. Since it is normally
generated, we know that $b-ae \geq 2g+1$. So Theorem
\ref{thm:minimalsection} guarantees that $L_C$ is normally
presented and hence $b-ae \geq 2g+2$.\\
$(3)$ By $(1)$ we may assume that $b-ae \geq 2g+1$. Now assume
that $L$ satisfies property $N_p$ for $p \geq 1$. We need to show
that $b-ae \geq 2g+1+p$. If $b-ae =2g+p$, then we can apply
Theorem \ref{thm:minimalsection} to $L$. Thus $L_C$ satisfies
property $N_p$ and hence $\mbox{deg}(L_C) \geq 2g+1+p$ which is a
contradiction.\\
$(4)$ If $a \geq 1$ and $b-ae \geq 2g+1+p$, then $L$ satisfies
property $N_p$ by Theorem \ref{thm:main}. The converse also holds
by $(3)$. \qed \\

\noindent As a corollary of Proposition
\ref{prop:normalgenerationsurface}, we also obtain the following:

\begin{corollary}
Let $X$ be a ruled surface over a curve $C$ of genus $g \geq 2$
with $e \geq 0$. For a very ample line bundle $L \in \mbox{Pic}X$
in the numerical class of $aC_0 +bf$, assume that $H^1 (X,L)=0$.
If $L$ is normally generated, then $L_C := L|_{C_0}$ is normally
generated.
\end{corollary}

\begin{proof}
From the short exact sequence $0 \rightarrow \mathcal{O}_X (-C_0)
\rightarrow \mathcal{O}_X \rightarrow \mathcal{O}_C \rightarrow
0$,\\ $H^1 (C,L_C)=0$. Thus $L_C$ is a nonspecial very ample line
bundle on $C$. In particular, $\mbox{deg}(L_C) = b-ae \geq g+3$.
Therefore it is easily checked that
\begin{equation*}
H^1 (X,L^{j+1} \otimes \mathcal{O}_X (-C_0))=0~~~\mbox{for
all}~~~~j \geq 1.
\end{equation*}
Now the assertion comes from the following Proposition
\ref{prop:normalgenerationsurface}.
\end{proof}

\noindent {\bf Remark 5.1.} When $X$ is a ruled surface over a
curve $C$ of genus $2$ with $e=-2$, $L \equiv C_0 +3f$ is very
ample, $H^1 (X,L^j)=0$ for all $j \geq 1$ and the complete linear
system $|L|$ defines an embedding $X \hookrightarrow \P^5$. Also
note that $L|_{C_0}$ is normally generated since
$\mbox{deg}(L|_{C_0})=L \cdot C_0 =5$. But $L$ is not normally
generated since $2$-normality fails to hold. For details, see
\cite{ABB}. Therefore $C_0$ is not the extremal curve with respect
to property $N_0$ and $L$. Thus it seems natural to investigate
extremal curves on ruled surfaces with negative $e$ as in the case
of elliptic ruled surfaces with
$e=-1$. \qed \\

\end{document}